\documentclass[centertags,leqno]{article}

\usepackage{amsmath}
\usepackage{amssymb}
\usepackage{latexsym}
\usepackage{amsxtra}
\usepackage{amscd}
\usepackage{theorem}
\usepackage{xypic}

\theoremstyle{plain}

{\theorembodyfont{\slshape}

        \newtheorem{thm}{Theorem}[section]
        \newtheorem{cor}[thm]{Corollary}
        \newtheorem{lem}[thm]{Lemma}
        \newtheorem{prop}[thm]{Proposition}

}

{\theorembodyfont{\rmfamily}

        \newtheorem{defn}[thm]{Definition}
        
        \newtheorem{rem}[thm]{Remark}

        \newtheorem{exa}[thm]{Example}

}

\renewcommand{\em}{\sl}

\newcommand{\proof}{{\bf Proof:\ }}

\newcommand{\Endproof}{\hspace*{\fill} $\Box$ \vspace{1ex} \noindent }

\makeatletter
\renewcommand{\subsection}{\@startsection{subsection}{2}%
        {\z@}{-3.25ex plus -1ex minus-.2ex}{-1em}{\bf}}
\makeatother

\newcommand{\PP}{\mathbb{P}}
\newcommand{\ZZ}{\mathbb{Z}}

\newcommand{\FF}{\mathbb{F}}

\newcommand{\X}{\mathcal{X}}
\newcommand{\Y}{\mathcal{Y}}

\newcommand{\Hom}{\mathop{\rm Hom}\nolimits}

\newcommand{\Spec}{\mathop{\rm Spec}\nolimits}

\title{Construction of covers in positive characteristic via degeneration}
\author{Irene I.\ Bouw}
\date{}
\begin{document}
\maketitle 
\begin{abstract} In this note we construct examples of covers of the 
projective line 
in positive characteristic such that every specialization is
inseparable. The result illustrates that it is not possible to
construct all covers of the generic $r$-pointed curve of genus zero
inductively from covers with a smaller number of branch points.
\\[2ex] 2000 Mathematical Subject Classification: Primary 14H30, 14H10
\end{abstract}

Let $k$ be an algebraically closed field of characteristic $p>0$. Let
$X=\PP^1_k$ and $G$ a finite group. We fix $r\geq 3$ distinct points
${\mathbf x}=(x_1, x_2, \ldots, x_r)$ on $X$. We ask whether there
exists a tame Galois cover $f:Y\to X$ with Galois group $G$ which is
branched at the $x_i$. If $p$ does not divide the order of $G$, then
the answer is well known. Namely, such a cover exists if and only if
$G$ may be generated by $r-1$ elements of order prime to $p$.

Suppose that $p$ divides the order of $G$. Then the existence of a
$G$-cover as above, depends on the position of the branch points
$x_i$. (See, for example, \cite[Lemma 6]{An}.) In this note we
restrict to the case that $(X; {\mathbf x})$ is the generic $r$-pointed
curve of genus zero.  A more precise version of the existence question
in positive characteristic is whether there exists a $G$-Galois cover
of $(X; {\mathbf x})$ with given ramification type (see for example
\cite{An}). For the particular kinds of groups we consider here, we
define the ramification type in \S \ref{metasec}.

Osserman (\cite{Osserman}) proves (non)existence of covers in positive
characteristic, for certain ramification types.  His method is roughly
as follows. First, he proved results for covers branched at $r=3$
points. In this case his results are strongest. Using the case $r=3$,
he then constructs {\sl admissible covers} of degenerate curves which
deform to covers of smooth curves (see \S \ref{admsec} for a
definition).

Suppose we are given a tame $G$-Galois cover $\pi$ of $(X=\PP^1_k; {\mathbf
 x})$.  Osserman asks (\cite[\S 6]{Osserman}) whether there exists a
 degeneration $(\bar{X}, \bar{{\mathbf x}})$ of $(X; {\mathbf x})$
 such that $\pi$ specializes to an admissible cover of $(\bar{X},
 \bar{{\mathbf x}})$. If such a degeneration exists, he says that
 $\pi$ has {\sl a good degeneration}. Covers which admit a good
 degeneration are exactly those which may be shown to exist
 inductively from the existence of covers with less branch points.  To
 goal of this note is to produce covers which do not have a good
 degeneration. We show that such covers exist with arbitrary large
 number of branch points.

\section{Meta-abelian covers}\label{metasec}
In this section, we recall a result from \cite{cyclic} on the
existence of tame Galois covers with Galois group $G\simeq
(\ZZ/p)^n\rtimes \ZZ/m$.  Let $p\neq 2$ be a prime and $m$ be an
integer prime to $p$. Let $f$ be the order of $p \pmod{m}$. We suppose
that $k$ is an algebraically closed field of characteristic
$p$.

Let ${\mathbf x}=(x_1, \ldots x_r)$ be $r$ distinct $k$-rational
points of $X=\PP^1_{k}$. Let ${\mathbf a}=(a_1, \ldots, a_r)$ be an
$r$-tuple of integers with $0< a_i<m$ and $\sum a_i\equiv 0\bmod{m}.$
Suppose moreover that $\gcd(m, a_1, \ldots, a_r)=1$.  Let $g:Z\to X$
be the $m$-cyclic cover of {\sl type} $({\mathbf x};{\mathbf a})$
(\cite{cyclic}), i.e.\ $Z$ is the complete nonsingular curve given by
the equation
\[
z^m=\prod_{i:x_i\neq \infty} (x-x_i)^{a_i}
\]
and $g:(x, z)\mapsto x$.  We denote by $\sigma(Z)$ the $p$-rank of
$Z$. Then $\sigma(Z)=\dim_{\FF_p} V$, where $V:=\Hom(\pi_1(Z),
\ZZ/p)$. Since $\ZZ/m\ZZ$ acts on $V$, there exists a tame $G:=V\rtimes
\ZZ/m\ZZ$-Galois cover $\pi:Y\to X$ which factors through $Z$.

The following proposition gives an upper bound on $\sigma(Z)$ which
is attained if the branch points $x_i$ are sufficiently general. For
 a more precise version, we refer to
\cite{cyclic}. See also \cite{hypergeo} for the case $f=1$.  For every
integer $a$, we denote by $\langle a\rangle$ the unique integer with
$\langle a\rangle\equiv a\pmod{m}$ and $0< \langle a\rangle<m$.

Let $\chi:\ZZ/m\ZZ\to\FF_{p^f}^\times$ be a nontrivial, irreducible
character. Let $I=\{1, \ldots, m-1\}/\sim$, where $i\sim p^j i$. Then
$I$ corresponds to the set of nontrivial, irreducible characters 
$\ZZ/m\ZZ\to \FF_p^\times$. For every $i\in I$, we let
$n_i:=(p^f-1)/\gcd(i,m)$ be the number of elements of the equivalence
class of $i$.

\begin{prop}\label{HWprop}
\begin{itemize}
\item[(a)] We have that
\[
\sigma(Z)\leq B({\mathbf a}):=\sum_{i\in I}n_i\min_{0\leq i\leq f-1}
(r-1-\frac{1}{m}\sum_{j=1}^r \langle p^ia_j\rangle).
\]
\item[(b)]
Suppose that $p\geq m(r-3)$. There exists $x_1, \ldots, x_r\in
X=\PP^1_k$ such that
\[
\sigma(Z)= B({\mathbf a}).
\]
\end{itemize}
\end{prop}

\proof Part (a) is proved in \cite[Lemma 4.3]{cyclic}. Part (b)
follows from \cite[Theorem 6.1]{cyclic}.
\Endproof

In \cite{cyclic} one finds some variants of this result: under certain
additional hypotheses on the type, we may weaken the condition on $p$.

The number $r-1-(\sum_{j=1}^r \langle ia_j\rangle)/m$ is the dimension
of the $\chi^{-i}$th-eigenspace of $H^1(C, {\mathcal O}_C)$
(\cite{cyclic}). It is well-known that this number is an upperbound
for the dimension of the $\chi^{i}$th eigenspace of
$V\otimes_{\FF_p}\FF_{p^f}$. The following statement immediately follows
from Proposition \ref{HWprop}.

\begin{cor}\label{ordinarycor}
Let $g:Z\to X$ be an $m$-cyclic cover of type $({\mathbf x};{\mathbf
a})$, where $(X; {\mathbf x})$ is generic. Suppose that $p\geq m(r-3)$.
Define 
\[
\gamma(s)=\frac{1}{m}\sum_{t=1}^r
\langle s a_t\rangle.
\]
Then $Z$ is ordinary if and only $\gamma(s)=\gamma(p^is)$, for all $i$.
\end{cor}

\section{Degeneration}\label{admsec}
Let $R=k[[t]]$ be a discrete valuation ring of equal characteristic
$p$ and let $\X\to \Spec(R)$ be a semistable curve over $R$ whose
generic fiber is smooth. Let $x_1, \ldots, x_r:\Spec(R)\to \X$ be
disjoint section, which avoid the singularities of $X_0:=\X\times_R
k$.

\begin{defn}\label{gooddegendef} Let $\pi_K:Y_K\to X_K$ be a tame cover of 
smooth projective curves.  We say that $\pi_K$ has a {\sl good
 degeneration} if there exists a discrete valuation ring $R$ with
 fraction field $K$ and a finite morphism $\pi:\Y\to \X$ of semistable
 curves over $\Spec(R)$ with generic fiber $\pi_K$ such that the
 branch locus is \'etale over $\Spec(R)$ and the special fiber is
 separable. If this holds, we call $\pi_R:\Y\to\X$ (or also its
 special fiber $\pi_0:=\pi_R\otimes_R k:Y_0\to X_0$) a {\sl good
 degeneration} of $\pi_K$.
\end{defn}

Let $\pi_K:Y_K\to X_K$ be a tame cover of smooth projective curves which
 has a good degeneration. Let $\pi_R:\Y\to\X$ be as in the statement of
 Definition \ref{gooddegendef}. Then the special fiber
 $\pi_0:=\pi\otimes_R k:Y_0\to X_0$ is an {\sl admissible cover}. We
 recall the definition and refer to \cite[\S 2.1]{An} for a short
 introduction to admissible covers.    Let $\tau$ be any singularity
 of $Y_0$, and let $Y_1, Y_2$ be the (not necessarily different)
 irreducible components of $Y_0$ which intersect in $\tau$. Then we
 require that the canonical generators $h_i$ (with respect to some
 chosen system of roots of unity) of the stabilizer of $\tau\in Y_i$
 satisfy $h_1\cdot h_2=1$. (Recall that $h$ is a {\sl canonical
 generator} if there exists a local parameter $u$ of $\tau$ such that
 $h^\ast u=\zeta_n\cdot u$, where $n$ is the order of the stabilizer
 of $\tau$.)

Let $G=(\ZZ/p\ZZ)^n\rtimes\ZZ/m\ZZ$ and $R=k[[t]]$ and
$K=k((t))$. Suppose that $\pi_K:Y_K\to X_K$ is a tame $G$-Galois cover
which has a good degeneration. Let $\pi:\Y\to\X$ be a finite morphism as
in Definition \ref{gooddegendef}.  It is easy to see if $\pi_K$ has a
good degeneration, then there exists a good degeneration $\pi:\Y\to\X$
such that the special fiber $X_0$ of $\X$ consists of two projective
curves meeting in one point $\tau$. We denote these components by
$X_1$ and $X_2$. Write $S_i\subset \{1, \ldots, r\}$ for the indices
$j$ such that $x_j$ specializes to $X_i$. We write $\pi_i:Y_i\to X_i \
(i=1,2)$ for the restriction of $\pi_0:=\pi\otimes_R k$ to $X_i$.

Let $Z_K=Y_K/(\ZZ/p\ZZ)^n$ and write $g_K:Z_K\to X_K$ for the
$m$-cyclic cover associated to $\pi_K$. Let $({\mathbf x}; {\mathbf a})$ be the
type of $g_K$.

\begin{lem}\label{admlem} Let $\pi_K:Y_K\to X_K$ be a $G$-Galois cover as above. 
Suppose that $X_0:=\X\otimes_R k$ consists of two irreducible
components $X_1, X_2$, as above. Then $\pi_i:Y_i\to X_i \ (i=1,2)$ has
type
\[
(( x_j)_{j\in S_i}\cup (\tau)\, ;\, (a_j)_{j\in S_j}\cup
(\sum_{j\not\in S_j} a_j)).
\]
\end{lem}

\proof
This follows immediately from the definition of the type.
\Endproof

A well-known result of formal patching (\cite{HS}, \cite{admissible})
states that every tame admissible cover may be deformed to a cover of
smooth curves. This may be used to produce examples of covers which
have a good degeneration. For example, one easily checks the
following. Let $\pi:Y\to X$ be a $G$-Galois cover of type $({\mathbf x};
{\mathbf a})$, where $(X; {\mathbf x})$ is the generic $r$-pointed
curve of genus zero. If there exists $1\leq i<j\leq r$ such that
$a_i+a_j=m$, then $\pi$ has a good degeneration.  We give an easy
example of a cover which does not have a good degeneration.

\begin{exa}\label{m=5exa} Let $m=5$ and let $p\equiv -1\pmod{m}$. Then the
 order, $f$, of $p$ in $\ZZ/m\ZZ^\ast$ is $2$. We consider ${\mathbf
a}=(1,1,1,2)$. One computes that $B({\mathbf a})=2$.  Proposition
\ref{HWprop} implies that for $p$ sufficiently large there exists a
tame $G=(\ZZ/p\ZZ)^2\rtimes \ZZ/m\ZZ$-Galois cover $\pi:Y\to \PP^1$
branched at $4$ points which factors through a cover of type ${\mathbf
a}$.  In fact, \cite[Proposition 7.8]{cyclic} implies that we do not
need the lower bound on $p$ in this case.

Let $\pi_0:Y_0\to X_0$ be a degeneration of $\pi$. Then $X_0$ consists
of $2$ irreducible components, which we denote by $X_1$ and $X_2$. To
each of these components specialize two of the points $x_1, \ldots,
x_4$. Lemma \ref{admlem} implies that (up to renumbering) the
restrictions $\pi_1$ and $\pi_2$ of $\pi_0$ would have type ${\mathbf
a}_1=(1,1,3)$ and ${\mathbf a}_2=(2,1,2)$. One computes that
$B({\mathbf a}_i)=0$ for $i=1,2$. Hence $\pi_0$ is
inseparable.Therefore $\pi$ does not have a good degeneration.
\end{exa}

\begin{rem}\label{f=12rem} Suppose that $p\equiv 1\pmod{m}$ and let 
$(X=\PP^1_k; {\mathbf x})$ be the generic $r$-pointed curve of genus
zero.  Then it is shown in \cite[Proposition 7.4]{cyclic} that every
$m$-cyclic cover of $(X; {\mathbf x})$ has a good
degeneration. Moreover, in this case we have that $B({\mathbf
a})=g(Z)$, for every type $m$-cyclic cover $Z\to X$ of type $({\mathbf
x};{\mathbf a})$.

In the case that $p\equiv -1\pmod{m}$ it is shown in \cite[Proposition
  7.8]{cyclic} that every $m$-cyclic cover of $(X=\PP^1_k; {\mathbf
  x})$ has a good degeneration, provided that the number of branch
points is at least $5$. The proof of this result relies essentially on
the fact that the group scheme $J(Z)[p]$ of $p$-torsion points of an
$m$-cyclic cover of $(X=\PP^1_k; {\mathbf x})$ is self-dual under
Cartier duality. The examples from \S \ref{degensec} suggest that a
similar result does not hold for $f>2$, see Remark \ref{flargerem}.
\end{rem}

\section{Covers without a good degeneration}\label{degensec}
In this section, we produce examples of Galois covers which do not
have a good degeneration. Let $f$ be an odd prime and put
$\alpha:=2$. Define $m:=\alpha^f-1=1+\alpha+\cdots+\alpha^{f-1}$.
 We define
\[
{\mathbf a}=(1, \alpha, \alpha^2, \ldots, \alpha^{f-1}).
\]
 We suppose that $(X=\PP^1_k; {\mathbf x})$ is the generic
$f$-pointed curve of genus zero and let $g:Z\to X=\PP^1_k$ be the
$m$-cyclic cover of type $({\mathbf x}; {\mathbf a})$.

 As in \S \ref{metasec}, we define
$\gamma(s)=(\sum_{t=0}^{f-1}\langle s\alpha^t\rangle)/m.$

\begin{lem}\label{combilem}
Let $S\subsetneq \{0, \ldots, f-1\}$ and $s:=\sum_{j\in S} \alpha^j$. Then
\[
\gamma(s)=|S|, \qquad \gamma(m-s)=f-|S|.
\]
\end{lem}

\proof
Let $s$ be as in the statement of the lemma.
The definition of $m$ implies that $\alpha^f\equiv 1\pmod{m}$. Therefore 
\[
\langle s \alpha^i\rangle=\sum_{j\in S} \alpha^{i+j},
\]
where the powers of $\alpha$ should be read modulo $f$. This implies that
\[
\gamma(s)=\frac{1}{m}\sum_{t=0}^{f-1}\sum_{j\in S}\langle
\alpha^{j+t}\rangle=\frac{1}{m}\sum_{j\in S}(1+\alpha+\cdots
+\alpha^{f-1})=|S|.
\] 
The second statement follows immediately from the first statement and
the definition of $m$.  \Endproof

\begin{lem}\label{boundlem}
Let $p\geq m(f-3)$ be a prime such that $p^f\equiv 1\pmod{m}$.
Then $Z$ is ordinary, i.e.\
\[
b:=B({\mathbf a})=g(Z)=(f-1)(m-1)/2.
\]
In particular, there exists a tame $G:=(\ZZ/p\ZZ)^b\rtimes
\ZZ/m\ZZ$-Galois cover  $\pi:Y\to\PP^1_k$ of type ${\mathbf a}$.
\end{lem}

\proof The assumption on $p$ implies that $p\equiv \alpha^i\pmod{m}$,
for some $i$. Therefore $\gamma(sp^i)=\gamma(s)$, for all $i$. The
statement now follows immediately from Corollary \ref{ordinarycor}.
\Endproof

We now suppose that the order of $p$ in $\ZZ/m\ZZ^\ast $ is $f$.  The
goal of this section is to show that the cover $\pi$ from Lemma
\ref{boundlem} does not have a good degeneration.  It suffices to show
that every degeneration $g_0:Z_0\to X_0$ of $g:Z\to X$ is
nonordinary. As remarked in \S \ref{admsec}, it suffices to consider
degenerations $g_0:Z_0\to X_0$ of $g:Z\to X$ such that fiber $X_0$
consists of two irreducible components $X_1$ and $X_2$ intersecting in
one point $\tau$.

Consider such a degeneration $g_0:Z_0\to X_0$.  We let $S_i\subset
\{1, \ldots, f\}$ be the subset of indices $j$ such that $x_j$
specializes to $X_i$ ($i=1,2$). We may assume that $2\leq |S_i|\leq
f-2$.

\begin{prop}\label{baddegprop}
Let $g_0:Z_0\to X_0$ be a degeneration of $g:Z\to X$. Then $Z_0$
is nonordinary.
\end{prop}

\proof We assume that $X_0$ consists of two irreducible components
$X_1$ and $X_2$ which intersect in one point $\tau$, and let $S_i$ be
as above. We write $g_i:Z_i\to X_i$ for the restriction of $g_0$ to $X_i$.
The curve $Z_0$ is ordinary if and only if $Z_i$ is ordinary, for $i=1,2$.
 
We define
\[
\gamma_1(i)=\frac{1}{m}\left(\langle\sum_{j\not\in S_1}
i\alpha^j\rangle+ \sum_{j\in S_1} \langle\alpha^j\rangle).\right)
\]
These are the terms occuring in the bound for the cover $g_1:Z_1\to
X_1$ (Lemma \ref{admlem}). Note that $g_1$ is branched at $r_1+1$
points, namely the specialization of $x_j$ for $j\in S_1$ and the
singular point $\tau$.  It follows from Corollary \ref{ordinarycor}
that $Z_1$ is ordinary if and only if $\gamma_1(i)=\gamma_1(p^ji),$
for all $j$.

Let $s=\sum_{j\in S_1} \alpha^j$.  Put
\[
d_s=\gamma_1(s)-\sum_{j\in S_1}\langle s\alpha^j\rangle=\langle s
\sum_{j\not\in S_1} \alpha^j\rangle=\langle s(m-s)\rangle=\langle m-s^2\rangle.
\]
Then Lemma \ref{combilem} implies that
\[
\gamma(d_s)=\gamma(m-s^2)=f-\gamma(s^2).
\]
Since $f$ is odd, it follows that $\gamma(s^2)=\gamma(s)$. We conclude that
\[
\gamma(d_s)=f-\gamma(s)=f-|S_1|.
\]

We claim that there exists an $i$
such that $\gamma_1(sp^i)\neq \gamma_1(s)$. Namely,
\[
\sum_{i=0}^{f-1} \gamma_1(s\alpha^i)=|S_1|^2+f-|S_1|\equiv |S_1|(|S_1|-1)\pmod{f}.
\]
Since $2\leq |S_1|\leq f-2$, we conclude therefore that
$\sum_{i=0}^{f-1} \gamma_1(s\alpha^i)\not\equiv 0\pmod{f}$. This shows
that there exists an $i$ such that $\gamma_1(p^is)\neq \gamma_1(s)$.
Hence $Z_1$ is not ordinary.
\Endproof

As already remarked above, the following corollary immediately follows
from Proposition \ref{baddegprop}.

\begin{cor} Let $\pi:Y\to \PP^1$ be as in Corollary \ref{ordinarycor}. Then $\pi$ does not
 have a good degeneration.
\end{cor}

\begin{rem}\label{flargerem} It seems that Proposition \ref{baddegprop} is a 
special case of a much more general statement. Let $m$ be an odd
integer and let $\alpha$ be an element of order $n|(m-1)/2$ in
$\ZZ/m\ZZ^\ast$. Let $p$ be sufficiently large such that $p$ has order
$f|n$ in $\ZZ/m\ZZ^\ast$. Let $g:Z\to X$ be an $m$-cyclic cover of
type $({\mathbf x};{\mathbf a})$, where $(X; {\mathbf x})$ is the
generic $n$-pointed curve of genus $0$ and ${\mathbf a}=(1, \alpha,
\ldots, \alpha^{n-1})$. As in the proof of Corollary
\ref{ordinarycor}, one checks that $Z$ is ordinary. Explicit
computations suggest that the tame
$(\ZZ/p\ZZ)^{g(Z)}\rtimes\ZZ/m\ZZ$-Galois cover corresponding to $g$
does not have a good degeneration if $f$ is odd and strictly larger
than $1$. This suggests that for $f\geq 3$ there is no generalization
of the statement of Remark \ref{f=12rem}.
\end{rem}


\begin{thebibliography}{10}
\bibitem{cyclic} I.~I. Bouw.  \newblock The $p$-rank of ramified
covers of curves.  \newblock {\em Compositio Math.},126:295--322,
2001.

\bibitem{hypergeo}
I.~I. Bouw.
\newblock Reduction of the {H}urwitz space of metacyclic covers.
\newblock {\em Duke Math.\ J.}, 121:75--111, 2004.

\bibitem{HS}
D. Harbater and K. Stevenson.
\newblock {\em  J. Algebra}, 212:272--304, 1999.  


\bibitem{Osserman}
B. Osserman.
\newblock Linear series and existence of branched covers.
\newblock To appear in {\em Compositio Math.}



\bibitem{admissible}
S. Wewers.
\newblock Deformation of tame admissible covers of curves.
 \newblock In: {\em Aspects of Galois theory (Gainesville, FL, 1996)},
   London Math. Soc. Lecture Note Ser., 256: 239--282, 1999.
    

\bibitem{An}
 S. Wewers and I.~I. Bouw.
\newblock Alternating groups as monodromy groups in positive characteristic.
\newblock {\em Pacific J. Math.} 222:185--200, 2005.


\end{thebibliography}
\end{document}